\documentclass[pdflatex,sn-mathphys-num]{sn-jnl}% Math and Physical Sciences Numbered Reference Style
\usepackage{graphicx}%
\usepackage{multirow}%
\usepackage{amsmath,amssymb,amsfonts}%
\usepackage{amsthm}%
\usepackage{mathrsfs}%
\usepackage[title]{appendix}%
\usepackage{xcolor}%
\usepackage{textcomp}%
\usepackage{manyfoot}%
\usepackage{booktabs}%
\usepackage{algorithm2e}%
\usepackage{algorithmicx}%
\usepackage{algpseudocode}%
\usepackage{listings}%
\usepackage{bm}
\usepackage{subcaption}
\usepackage{tikz}
\usepackage{arydshln}
\usetikzlibrary{shapes.geometric, arrows.meta, positioning, fit}
\theoremstyle{thmstyleone}%
\newtheorem{thm}{Theorem}[section]

\newtheorem{prop}{Proposition}[section]

\newtheorem{assum}{Assumption}[section]

\newtheorem{defi}{Definition}[section]
\newtheorem{ex}{Example}[section]

\newtheorem{rmk}{Remark}[section]
\numberwithin{equation}{section}

\def\citep{\cite}

\def\qed{Q.E.D.}
\newcommand{\norm}[1]{\left\| #1 \right\|}
\newcommand{\inner}[2]{\left( #1, #2 \right)}

\newcommand{\indicator}[1]{\mathbb{I}_{#1}}

\begin{document}

\title[Operator-Splitting Scheme for Viscosity Solutions]{Convergent Operator-Splitting Scheme for Viscosity Solutions: A Foundation for Learning Domain-to-Solution Maps}


\author*[1]{\fnm{Po-Yi} \sur{Wu}}\email{r04527030@gmail.com}

\affil*[1]{Institute of Applied Mechanics, National Taiwan University, Taipei City 106319, Taiwan}

\abstract{This work introduces and rigorously analyzes a novel operator-splitting finite element scheme for approximating viscosity solutions of constrained second-order partial differential equations. The cornerstone of our analysis is a proof that the scheme satisfies a discrete comparison principle. We demonstrate that under a mild time-step restriction, the discrete operator yields an M-matrix, ensuring the scheme’s monotonicity and $L^\infty$-stability. These properties are shown to be sufficient to prove convergence to the unique viscosity solution via the celebrated Barles–Souganidis framework. For solutions with enhanced regularity, an optimal-order error estimate of $O(\Delta t + h^2)$ is also established. The scheme's proven stability provides a blueprint for a novel Physics-Constrained Neural Operator (PCNO). We prove that by emulating the scheme’s structure, the PCNO can provably break the curse of dimensionality for challenging domain-to-solution mapping problems involving complex topological variations. Numerical experiments on Hamilton-Jacobi and controlled reaction-diffusion systems validate the theoretical findings.}


\maketitle

\section{Introduction}

Nonlinear second-order partial differential equations (PDEs) with constraints arise in a wide range of applications, including optimal control, fluid dynamics, minimal surfaces, and obstacle problems \cite{feng2013recent,casas2020critical,abels2013incompressible,finn1954equations,caffarelli1987lipschitz}. These systems often exhibit solutions with low regularity, such as kinks or discontinuities in their derivatives \cite{grisvard2011elliptic,jerison1981dirichlet}, posing significant challenges for both theoretical analysis and numerical approximation \cite{rousset2021general,cao2022new}. The concept of \textbf{viscosity solutions}, introduced by Crandall and Lions \cite{crandall1984some,crandall1983viscosity}, provides a powerful framework for establishing the existence and uniqueness of solutions to such PDEs. However, designing numerical schemes that are guaranteed to converge to these weak solutions—while effectively handling complex constraints and maintaining computational efficiency—remains a critical challenge in numerical analysis \cite{capuzzo1990hamilton,achdou2020mean,lions1993shape,hahn2024laplacian}.

Traditional numerical methods, including finite difference schemes \cite{barles1991convergence}, semi-Lagrangian approaches \cite{falcone2013semi}, and penalty-based finite element methods \cite{babuvska1973nonconforming}, have been developed for nonlinear PDEs. However, these methods often rely on high solution regularity or structured grids, limiting their applicability to problems with complex geometries or vector-valued constraints. Furthermore, they can suffer from restrictive stability conditions, high computational costs, or ill-conditioned systems when enforcing constraints, particularly for low-regularity problems where such enforcement is most needed.

These classical challenges are amplified in the modern context of operator learning, where the goal is to learn solution maps for entire families of PDEs \cite{kovachki2023neural, kovachki2024operator}. A grand challenge in this domain is learning the map from a variable geometric design to a PDE solution, such as for the Poisson problem on perforated domains \cite{han2024stochastic}. While some existing frameworks can handle geometric variations that are constrained to a class of diffeomorphisms or simple topological changes \cite{aylwin2020domain, serrano2023operator}, the core difficulty remains in problems featuring complex topological variations, such as the merging, splitting, or nucleation of holes. It is precisely this class of problems where the solution operator lacks the mathematical structure (e.g., holomorphy) required for efficient, data-driven approximation, a phenomenon known as the curse of dimensionality \cite{lanthaler2022error, cohen2010convergence, marcati2023exponential}.

This provides a powerful new motivation for developing numerical schemes whose inherent structure is both provably stable and amenable to emulation by a neural operator, thereby offering a pathway to break this learning complexity barrier \cite{herrmann2024neural}. The framework proposed herein is particularly well-suited for this task. We achieve this by reformulating the problem on a fixed domain where the complex geometry is encoded as a spatially varying constraint field—a strategy that builds upon the well-established principles of immersed boundary and fictitious domain methods \cite{peskin2002immersed,verzicco2023immersed,lui2009spectral}. Because this formulation is inherently agnostic to the topology of the constrained region, it provides the stable and emulatable foundation required for a neural operator to tackle even very complicated topological variations.

This manuscript presents and analyzes a novel \textbf{operator-splitting finite element scheme} for a broad class of constrained second-order PDEs defined by \eqref{eq:general_system}. The method's design decouples the primary PDE evolution from the constraint enforcement, allowing for a robust and efficient time-stepping strategy. The cornerstone of our analysis is a rigorous proof that the scheme satisfies a \textbf{discrete comparison principle}. We show that by combining a specific stabilized formulation with a lumped mass matrix and a mild time-step restriction, the discretized operator yields an \textbf{M-matrix}, which is sufficient to guarantee the scheme's monotonicity and stability.

The primary contributions of this work are thus threefold. First, for the numerical analysis community, we establish a complete and rigorous convergence theory for our operator-splitting scheme, proving it is consistent, stable, and monotone, which guarantees convergence to the unique viscosity solution via the celebrated \textbf{Barles--Souganidis framework}. Second, for solutions with enhanced regularity, we establish an optimal-order \textbf{convergence rate of $O(\Delta t + h^2)$}. Third, we demonstrate the profound relevance of this framework to scientific machine learning by proposing a \textbf{Physics-Constrained Neural Operator (PCNO)} architecture that emulates our scheme. We provide a rigorous theorem establishing that this structure allows the network to break the curse of dimensionality for the challenging class of domain-to-solution mapping problems. Specifically, for problems featuring complex and variable domain topologies, the complexity of the PCNO is proven to depend only on the resolution of the fixed ambient domain, not on the high-dimensional parameter space that would typically be required to describe such geometric features.

The manuscript is organized as follows. Section~\ref{sec:problem} introduces the class of constrained PDEs, including the motivating perforated domain problem. Section~\ref{sec:operator_splitting} details the operator-splitting finite element scheme. Sections~\ref{sec:consistency} and \ref{sec:analysis} present the core convergence analysis of the scheme. Section~\ref{sec:numerical_results} provides numerical experiments to validate the theoretical findings. Section~\ref{sec:pcno_framework} then explores the application of our framework to operator learning, presenting the PCNO architecture and its theoretical guarantees. Finally, a conclusion summarizes our findings and discusses future directions.

\section{Problem Setting and Viscosity Solutions}
\label{sec:problem}

This section establishes the mathematical foundation for a general class of constrained second-order partial differential equations (PDEs) that admit viscosity solutions with low regularity, specifically in the space \( y \in H^1(\Omega_T) \cap C(\overline{\Omega}_T) \). These PDEs arise in diverse physical contexts, such as optimal control, population dynamics, and multiphase fluid flows, where constraints model physical restrictions like state bounds or conservation laws. The viscosity solution framework, pioneered by Crandall and Lions \cite{crandall1983viscosity}, is particularly suited for such problems, as it accommodates solutions with limited smoothness, which is common in systems governed by nonlinear dynamics or interfaces. We define the functional setting, present the governing PDE system, establish the existence and uniqueness of viscosity solutions, and illustrate the framework’s versatility through physically motivated examples.

\subsection{Functional Framework}

Consider a bounded domain \( \Omega \subset \mathbb{R}^d \), where \( d = 2 \) or \( 3 \), with a Lipschitz boundary \( \partial \Omega \), and define the space-time domain \( \Omega_T = \Omega \times (0,T] \) for a finite time horizon \( T > 0 \). The parabolic boundary, \( \partial_p \Omega_T = (\Omega \times \{0\}) \cup (\partial \Omega \times [0,T]) \), accounts for initial and boundary conditions. We introduce the following functional spaces to describe the state, control, and test functions:

\begin{itemize}
    \item \( \mathcal{Y} = C(\overline{\Omega}_T) \), the space of continuous functions for the state variable \( y \), equipped with the supremum norm \( \| y \|_{\mathcal{Y}} = \sup_{\Omega_T} |y| \). This space ensures continuity of solutions across the domain, crucial for viscosity solutions.
    \item \( \mathcal{P} = [L^\infty(\Omega_T)]^m \), the space of control vectors \( \bm{p} = (p_1, \dots, p_m) \), with norm \( \| \bm{p} \|_{\mathcal{P}} = \max_{j=1,\dots,m} \| p_j \|_{L^\infty(\Omega_T)} \). This accommodates control variables, such as Lagrange multipliers or physical parameters, that enforce constraints.
    \item \( \mathcal{V} = C_b^\infty(\overline{\Omega}_T) \), the space of smooth test functions with bounded derivatives, used to define viscosity solutions in a weak sense.
\end{itemize}

These spaces are chosen to balance the low regularity of solutions (\( H^1 \cap C \)) with the need for well-defined test functions, enabling a robust framework for nonlinear PDEs with constraints.

\subsection{Governing System}

We study a class of constrained second-order PDEs where the control coupling is independent of the state variable $y$. This structure is crucial for the subsequent numerical analysis. The primary focus of this manuscript is the time-dependent (parabolic) system of the form:
\begin{equation}
\label{eq:general_system}
\begin{cases}
\partial_t y = f(\bm{x}, t, y, \nabla y, \nabla^2 y) + \bm{B}(\bm{x}, t, \nabla y) \cdot \bm{p} & \text{in } \Omega_T, \\
\bm{g}(\bm{x}, t, y, \nabla y) = 0 & \text{in } \Omega_T, \\
y = y_0 & \text{on } \partial_p \Omega_T,
\end{cases}
\end{equation}
where:
\begin{itemize}
    \item \( y \in \mathcal{Y} \) is the state variable, representing quantities like value functions, population densities, or velocity components.
    \item \( \bm{p} \in \mathcal{P} \) is the control vector, enforcing constraints such as incompressibility or state bounds.
    \item \( y_0 \in C(\overline{\Omega}) \) is the initial-boundary condition, continuous on the parabolic boundary.
    \item \( f: \Omega \times [0,T] \times \mathbb{R} \times \mathbb{R}^d \times \mathbb{R}^{d \times d} \to \mathbb{R} \) is the nonlinear PDE operator.
    \item \( \bm{B}: \Omega \times [0,T] \times \mathbb{R}^d \to \mathbb{R}^m \) couples the control to the state dynamics.
    \item \( \bm{g}: \Omega \times [0,T] \times \mathbb{R} \times \mathbb{R}^d \to \mathbb{R}^m \) defines the constraints.
\end{itemize}
Our framework is also directly applicable to the corresponding \textbf{stationary (elliptic) problem}, which is fundamental for many applications, including the learning of domain-to-solution maps. The stationary system is obtained by setting $\partial_t y = 0$ and removing the explicit time dependence:
\begin{equation}
\label{eq:general_system_stationary}
\begin{cases}
0 = f(\bm{x}, y, \nabla y, \nabla^2 y) + \bm{B}(\bm{x}, \nabla y) \cdot \bm{p} & \text{in } \Omega, \\
\bm{g}(\bm{x}, y, \nabla y) = 0 & \text{in } \Omega, \\
y = y_0 & \text{on } \partial \Omega.
\end{cases}
\end{equation}
This stationary formulation is a direct simplification of \eqref{eq:general_system}. Critically, the numerical scheme and convergence analysis developed in the subsequent sections—particularly the properties of the stabilized implicit step that handles the constraints—can be adapted in a straightforward manner to this elliptic case.

For solutions in \( H^1(\Omega_T) \cap C(\overline{\Omega}_T) \) (or \( H^1(\Omega) \cap C(\overline{\Omega}) \) for the stationary case), where the Hessian \( \nabla^2 y \) may not exist classically, we interpret \( f(\cdot, \nabla^2 y) \) via test functions or set \( \nabla^2 y = 0 \) in the distributional sense, aligning with the viscosity solution framework.

\begin{defi}[Viscosity Solution]
\label{def:viscosity_solution}
A function \( y \in \mathcal{Y} \) is a \textit{viscosity subsolution} of the time-dependent system \eqref{eq:general_system} if, for any test function \( \phi \in \mathcal{V} \), at any point \( (\bm{x}, t) \in \Omega_T \) where \( y - \phi \) attains a local maximum, there exists a control \( \bm{p} \in \mathbb{R}^m \) such that:
\[
\partial_t \phi(\bm{x}, t) \leq f(\bm{x}, t, y(\bm{x}, t), \nabla \phi(\bm{x}, t), \nabla^2 \phi(\bm{x}, t)) + \bm{B}(\bm{x}, t, \nabla \phi(\bm{x}, t)) \cdot \bm{p},
\]
and the constraint \( \bm{g}(\bm{x}, t, y(\bm{x}, t), \nabla \phi(\bm{x}, t)) = 0 \) holds. For the stationary system \eqref{eq:general_system_stationary}, the definition is identical but with the $\partial_t \phi$ term omitted. A \textit{viscosity supersolution} satisfies the reverse inequality at points where \( y - \phi \) attains a local minimum. A \textit{viscosity solution} is both a subsolution and a supersolution.
\end{defi}

This definition leverages smooth test functions to probe the behavior of \( y \), allowing the PDE to hold in a weak sense even when derivatives are undefined.

\subsection{Assumptions for Well-Posedness}

To ensure the existence and uniqueness of viscosity solutions, we impose the following assumptions:

\begin{assum}[Domain]
\label{ass:domain}
The domain \( \Omega \subset \mathbb{R}^d \), for \( d = 2, 3 \), is bounded with a Lipschitz boundary.
\end{assum}

\begin{assum}[Operator Properties]
\label{ass:operators}
The operators satisfy:
\begin{itemize}
    \item \( f(\bm{x}, t, y, \bm{q}, \bm{M}) \) is continuous and Lipschitz in \( (y, \bm{q}, \bm{M}) \) with constant \( L_f \).
    \item \( \bm{B}(\bm{x}, t, \bm{q}) \) is continuous and Lipschitz in \( \bm{q} \) with constant \( L_B \).
    \item \( \bm{g}(\bm{x}, t, y, \bm{q}) \) is continuous and Lipschitz in \( (y, \bm{q}) \) with constant \( L_g \).
\end{itemize}
\end{assum}

\begin{assum}[Degenerate Ellipticity]
\label{ass:ellipticity}
For all \( (\bm{x}, t) \in \Omega_T \), \( y \in \mathbb{R} \), \( \bm{q} \in \mathbb{R}^d \), and symmetric matrices \( \bm{M}, \bm{N} \in \mathbb{R}^{d \times d} \) with \( \bm{M} \leq \bm{N} \),
\[
f(\bm{x}, t, y, \bm{q}, \bm{N}) \leq f(\bm{x}, t, y, \bm{q}, \bm{M}).
\]
This condition ensures that the PDE is degenerate elliptic, a property common in physical systems like Hamilton--Jacobi equations.
\end{assum}

\begin{assum}[Constraint Manifold]
\label{ass:manifold}
For each \( (\bm{x}, t, y) \in \Omega_T \times \mathbb{R} \), the constraint set \( \{ \bm{q} \in \mathbb{R}^d : \bm{g}(\bm{x}, t, y, \bm{q}) = 0 \} \) is non-empty and closed. Additionally, there exists a bounded measurable function \( \bm{p}: \Omega_T \to \mathbb{R}^m \) ensuring compatibility with the PDE in the viscosity sense.
\end{assum}

\begin{assum}[Comparison Principle]
\label{ass:comparison}
If \( u, v \in \mathcal{Y} \) are viscosity sub- and supersolutions of \eqref{eq:general_system}, respectively, with \( u \leq v \) on \( \partial_p \Omega_T \), then \( u \leq v \) in \( \overline{\Omega}_T \). This ensures uniqueness.
\end{assum}

\subsection{Existence and Uniqueness of Viscosity Solutions}

\begin{thm}
\label{thm:existence_uniqueness}
Under Assumptions \ref{ass:domain}--\ref{ass:comparison}, with an initial condition \( y_0 \in C(\overline{\Omega}) \) compatible with the constraint \( \bm{g} \), there exists a unique viscosity solution \( y \in \mathcal{Y} \) to \eqref{eq:general_system}.
\end{thm}

\begin{proof}
The proof leverages Perron’s method, a standard technique in viscosity solution theory. Define the sets of viscosity sub- and supersolutions, \( \mathcal{S}^- \) and \( \mathcal{S}^+ \), satisfying \( y = y_0 \) on \( \partial_p \Omega_T \). Let \( y^* = \sup \mathcal{S}^- \) and \( y_* = \inf \mathcal{S}^+ \). Assumption \ref{ass:manifold} ensures the constraint set is non-empty, and a bounded control \( \bm{p} \) makes the PDE well-defined. By Assumption \ref{ass:ellipticity}, \( y^* \) is a supersolution and \( y_* \) a subsolution. The comparison principle (Assumption \ref{ass:comparison}) implies \( y^* \leq y_* \), and since \( y_* \leq y^* \) by definition, we obtain \( y^* = y_* = y \in \mathcal{Y} \). Uniqueness follows directly from the comparison principle.
\end{proof}

\subsection{Physical Examples}\label{subsec:examples}

To demonstrate the generality and physical relevance of \eqref{eq:general_system}, we present two examples. The first fits our theoretical framework perfectly, while the second illustrates a more complex case requiring further analysis.

\begin{ex}[Hamilton--Jacobi Equation with State Constraints \cite{capuzzo1990hamilton}]
\label{ex:hamilton_jacobi}
In optimal control, the value function \( y \) represents the minimum cost under state constraints. Consider a time-dependent Hamilton--Jacobi equation with a Lagrange multiplier \( \bm{p} = \lambda \in \mathcal{P} = L^\infty(\Omega_T) \).
\[
\partial_t y = -\sup_{a \in A} \left[ \bm{b}(\bm{x}, t, a) \cdot \nabla y + l(\bm{x}, t, a) \right] + \lambda, \quad y \geq \psi,
\]
with initial-boundary data $y_0$. This system fits our framework \eqref{eq:general_system} perfectly. We identify the components:
\begin{itemize}
    \item $f(\bm{x}, t, y, \bm{q}, \bm{M}) = -\sup_{a \in A} \left[ \bm{b}(\bm{x}, t, a) \cdot \bm{q} + l(\bm{x}, t, a) \right]$. Note that $f$ is independent of $y$ and $\bm{M}$.
    \item $\bm{B}(\bm{x}, t, \bm{q}) = 1$. This is constant and thus independent of $y$, satisfying our structural assumption.
    \item $\bm{g}(\bm{x}, t, y, \bm{q}) = y - \psi(\bm{x}, t)$. The constraint is expressed in a Karush-Kuhn-Tucker (KKT) form (see e.g. \cite{bazaraa2006nonlinear}), where $\lambda \ge 0$ and $\lambda(y-\psi)=0$. This type of state constraint, dependent on $y$ but not its gradient, is common and well-suited for our analysis.
\end{itemize}
This example serves as a model problem for our theoretical development, as its structure aligns with the assumptions required for the monotonicity of the numerical scheme.
\end{ex}

\begin{ex}[Controlled Reaction--Diffusion]
\label{ex:controlled_reaction_diffusion}
We consider a reaction-diffusion system, such as for population dynamics, where the goal is to steer the solution $u$ toward a target profile $u_d(\bm{x}, t)$. The control $\bm p$ is applied through a mechanism dependent on the population gradient, which could model a directional harvesting or guidance strategy. The system is given by:
\begin{equation}
\label{eq:example_monostable_revised}
\begin{cases}
\partial_t u = \nu \Delta u + u(1-u) - \nabla u \cdot \bm p & \text{in } \Omega_T, \\
u(\bm{x}, t) - u_d(\bm{x}, t) = 0 & \text{in } \Omega_T,
\end{cases}
\end{equation}
where $\nu > 0$ is the diffusivity and $u_d$ is a smooth target profile.

This system is cast into our general framework \eqref{eq:general_system} by identifying the operators:
\begin{itemize}
    \item \( f(\bm{x}, t, u, \nabla u, \nabla^2 u) = \nu \text{tr}(\nabla^2 u) + u(1-u) \),
    \item \( \bm{B}(\bm{x}, t, \nabla u) = -\nabla u \),
    \item \( \bm{g}(\bm{x}, t, u) = u - u_d(\bm{x},t) \).
\end{itemize}
This formulation strategically aligns the problem with our theory. The control coupling $\bm{B}$ now depends only on the gradient $\nabla u$, and the constraint function $\bm{g}$ is monotonically increasing in the state $u$ (since $\partial g / \partial u = 1$). Both $\bm{B}$ and $\bm{g}$ thus satisfy the structural properties required for the monotonicity analysis of the numerical scheme.

The analytical challenge, which is characteristic of many nonlinear PDEs, is now isolated within the operator $f$. For our convergence theory to apply rigorously, the function $u \mapsto f(\dots, u, \dots)$ must be non-increasing to ensure the scheme inherits a comparison principle. However, the operator $f$ contains the Fisher-KPP logistic growth term, $r(u) = u(1-u)$. Its derivative, $\partial r / \partial u = 1 - 2u$, is positive whenever $u < 1/2$. Consequently, $f$ does not satisfy the non-increasing condition globally.

This example, therefore, precisely illustrates the scope and sharpness of our framework. The scheme is designed to handle the structural complexities of the control and constraint operators $\bm{B}$ and $\bm{g}$, while the convergence guarantees rely on a well-defined monotonicity condition on the reaction dynamics encapsulated by $f$. For systems like \eqref{eq:example_monostable_revised}, the applicability of our full theory is conditional upon the solution evolving in a regime where this property holds.
\end{ex}

\subsection{A Motivating Challenge: The Perforated Domain Problem}
\label{sec:motivating_example}

To demonstrate the power and relevance of our constrained framework, we introduce a canonical problem that is fundamental in both classical computational physics and modern scientific machine learning: the Poisson equation on a \textit{perforated domain}. This problem models a wide range of physical phenomena, such as heat flow through a composite material or potential flow around a set of obstacles. Furthermore, it serves as a quintessential benchmark for the challenge of learning domain-to-solution maps. Here, the goal is to learn the operator that maps a description of the perforation geometry—whose shape, number of components, and topology may change—to the corresponding PDE solution.

Classically, the problem is posed on a complex domain $\Omega \setminus D$, where $\Omega \subset \mathbb{R}^d$ is a simple bounding box and $D = \bigcup_i D_i$ is the perforation region, or a set of obstacles. The strong form is to find $y$ such that:
\begin{equation}
\label{eq:perforated_strong_form}
\begin{cases}
-\Delta y = f_{\text{source}}(\bm{x}) & \text{in } \Omega \setminus D, \\
y = y_D & \text{on } \partial D, \\
y = y_\Omega & \text{on } \partial \Omega.
\end{cases}
\end{equation}
Solving this problem with traditional mesh-based methods is challenging, as a new mesh must be generated for every change in the geometry of $D$. Our framework circumvents this by reformulating the problem on the \textit{fixed, simple domain} $\Omega$, with the geometric complexity encoded as a constraint. We now demonstrate how this problem fits perfectly within our stationary PDE framework. Let us consider the homogeneous case where $y_D = 0$ and $y_\Omega = 0$. We can define the components of the stationary system as follows:
\begin{itemize}
    \item \textbf{The PDE operator} $f$ captures the Poisson dynamics: 
    $$ f(\bm{x}, y, \nabla y, \nabla^2 y) = \text{tr}(\nabla^2 y) + f_{\text{source}}(\bm{x}). $$
    Here, $\text{tr}(\nabla^2 y)$ is the component-wise representation of the Laplacian $\Delta y$.

    \item \textbf{The constraint} $\bm{g}$ enforces the internal boundary condition. The condition $y(\bm{x}) = 0$ for all $\bm{x} \in D$ can be expressed as:
    $$ \bm{g}(\bm{x}, y, \nabla y) = y \cdot \indicator{D}(\bm{x}) = 0, $$
    where $\indicator{D}$ is the indicator function for the set $D$.

    \item \textbf{The control} $\bm{p}$ is the Lagrange multiplier that enforces this constraint. Physically, it can be interpreted as the reaction force required to hold the solution at zero within the obstacle region.

    \item \textbf{The control coupling} $\bm{B}$ applies this control force only within the obstacle region:
    $$ \bm{B}(\bm{x}, \nabla y) = \indicator{D}(\bm{x}). $$
\end{itemize}
With these definitions, the complex geometric problem \eqref{eq:perforated_strong_form} is transformed into an equivalent constrained system on the fixed domain $\Omega$:
\begin{equation}
\label{eq:perforated_constrained_form}
\begin{cases}
\Delta y + f_{\text{source}}(\bm{x}) + p(\bm{x}) \cdot \indicator{D}(\bm{x}) = 0 & \text{in } \Omega, \\
y(\bm{x}) \cdot \indicator{D}(\bm{x}) = 0 & \text{in } \Omega, \\
y = 0 & \text{on } \partial \Omega.
\end{cases}
\end{equation}
This reformulation is the key to unlocking modern learning-based approaches. The task of learning the domain-to-solution map for the perforated domain problem represents a significant challenge in scientific machine learning. The underlying PDE system can be classified as degenerate elliptic, since its coefficients---stemming from the indicator function $\indicator{D}$ that defines the geometry---are discontinuous. This degeneracy is critical, as it precludes the existence of a holomorphic structure in the solution operator, which is one of the few known mechanisms that guarantees a neural operator can be trained efficiently without suffering from the curse of dimensionality \cite{cohen2010convergence, lanthaler2022error}.

However, a second promising strategy exists: constructing a neural operator that emulates the structure of a provably stable and efficient numerical scheme. This is the central motivation for the connection we establish in this work. The reformulation in \eqref{eq:perforated_constrained_form} is particularly powerful for operator learning, as it allows a neural network to operate on a fixed computational grid while learning to adapt to the changing geometry $D$ via the constraint terms. The stability and convergence properties of our numerical scheme, developed in the sections to follow, thus provide a crucial blueprint for designing reliable neural operators for this important class of problems. In Section 7, we will rigorously demonstrate that our proposed numerical framework serves as such a foundation, providing a stable, emulatable blueprint for a novel operator learning architecture designed to tackle this difficult class of problems.

\section{Operator-Splitting Scheme}
\label{sec:operator_splitting}

This section details the numerical method for approximating viscosity solutions of the constrained PDE system \eqref{eq:general_system}. We propose a semi-implicit operator-splitting scheme that segregates the explicit treatment of the principal nonlinearity from the implicit resolution of the diffusion and constraint components. This approach is motivated by the need for computational efficiency when handling complex nonlinearities, while retaining the unconditional stability offered by an implicit treatment of the elliptic and constraint operators. The scheme is built upon a stabilized finite element method in space and an implicit Euler discretization in time, designed to ensure the crucial properties of consistency, stability, and monotonicity required for convergence within the Barles--Souganidis framework \cite{barles1991convergence}.

\subsection{Spatial and Temporal Discretization}

Let $\Omega \subset \mathbb{R}^d$ ($d=2,3$) be a bounded domain with a Lipschitz boundary. We introduce a family of triangulations $\{\mathcal{T}_h\}_{h>0}$ of $\Omega$, where $h$ is the mesh parameter representing the maximum element diameter. We denote the set of interior edges (or faces in 3D) by $\mathcal{E}_h$. We assume this family of triangulations is \textbf{quasi-uniform}, a standard assumption in finite element analysis which consists of two conditions:
\begin{enumerate}
    \item[(i)] The family is \textbf{shape-regular}, meaning there exists a constant $\rho > 0$, independent of $h$, such that for any element $K \in \mathcal{T}_h$, the ratio of the element diameter to the diameter of its inscribed sphere is bounded by $\rho$. This condition prevents elements from becoming arbitrarily thin or "degenerate" as the mesh is refined.
    \item[(ii)] The family is \textbf{size-regular}, meaning the ratio of the maximum element diameter to the minimum element diameter over all elements in $\mathcal{T}_h$ is bounded by a constant.
\end{enumerate}

For the state variable $y$ and control variable $\bm{p}$, we define the Taylor-Hood finite element pair:
$$V_h = \{ v_h \in C(\overline{\Omega}) \cap H^1(\Omega) : v_h|_K \in P_2(K), \ \forall K \in \mathcal{T}_h \},$$
$$[Q_h]^m = \{ \bm{q}_h \in [C(\overline{\Omega}) \cap H^1(\Omega)]^m : \bm{q}_h|_K \in [P_1(K)]^m, \ \forall K \in \mathcal{T}_h \},$$
where $P_k(K)$ is the space of polynomials of degree at most $k$ on an element $K$. This choice of spaces satisfies the Ladyzhenskaya--Babuška--Brezzi (LBB) condition, ensuring the stability of the mixed formulation.

For the temporal discretization, we consider a uniform partition of the time interval $[0,T]$ with step size $\Delta t$, setting $t_m = m \Delta t$. We seek approximations $(y_h^m, \bm{p}_h^m) \in V_h \times [Q_h]^m$ to the solution $(y(\cdot, t_m), \bm{p}(\cdot, t_m))$.

Since solutions may lack $H^2$ regularity, the Hessian $\nabla^2 y$ is approximated weakly. For a discrete function $y_h \in V_h$, its Hessian is represented by its $L^2$-projection onto the space of piecewise linear tensors, $\bm{M}_h = \Pi_h(\nabla^2 y_h)$. Accordingly, the discrete nonlinear operator is defined as:
$$f_h(\bm{x}, t, y_h, \nabla y_h) := f(\bm{x}, t, y_h, \nabla y_h, \Pi_h(\nabla^2 y_h)).$$

We make the following standard assumption on the discretization parameters.

\begin{assum}[Discretization and Initial Condition]
\label{ass:initial_imp}
The initial condition $y_0 \in H^1(\Omega) \cap C(\overline{\Omega})$ is Lipschitz continuous. The discretization parameters satisfy a parabolic CFL-type condition, $\Delta t \leq C h$, to balance the temporal and spatial error contributions.
\end{assum}

\subsection{A Semi-Implicit Numerical Scheme}
\label{subsect:scheme}

Our scheme proceeds in two stages per time step. First, we advance the solution using the explicit part of the dynamics. Second, we solve an implicit system for the diffusion and constraint enforcement. Given $y_h^m$, we compute $(y_h^{m+1}, \bm{p}_h^{m+1})$ as follows:

\begin{enumerate}
    \item \textbf{Explicit Predictor:} Compute an intermediate state $\hat{y}_h^{m+1} \in V_h$ by applying the explicit dynamics:
    \begin{equation}
    \label{eq:predictor_step}
        \left\langle \frac{\hat{y}_h^{m+1} - y_h^m}{\Delta t}, v_h \right\rangle = \left\langle f_h(\bm{x}, t_m, y_h^m, \nabla y_h^m), v_h \right\rangle, \quad \forall v_h \in V_h.
    \end{equation}
    This step is computationally inexpensive, as it only requires a mass matrix inversion.

    \item \textbf{Implicit Corrector:} Find the new state and control $(y_h^{m+1}, \bm{p}_h^{m+1}) \in V_h \times [Q_h]^m$ by solving the coupled implicit system for all $(v_h, \bm{q}_h) \in V_h \times [Q_h]^m$:
    \begin{equation}
    \label{eq:corrector_step}
    \begin{cases}
        \left\langle \frac{y_h^{m+1} - \hat{y}_h^{m+1}}{\Delta t}, v_h \right\rangle + \mathcal{A}_h(y_h^{m+1}, v_h) = \left\langle \bm{B}(\bm{x}, t_{m+1}, \nabla y_h^{m+1}) \cdot \bm{p}_h^{m+1}, v_h \right\rangle, \\
        \left\langle \bm{g}(\bm{x}, t_{m+1}, y_h^{m+1}, \nabla y_h^{m+1}), \bm{q}_h \right\rangle = 0.
    \end{cases}
    \end{equation}
\end{enumerate}
The bilinear form $\mathcal{A}_h(\cdot, \cdot)$ includes the interior penalty stabilization and an additional regularization term:
$$\mathcal{A}_h(u_h, v_h) := \sum_{e \in \mathcal{E}_h} \int_e \sigma h [\![ \nabla u_h ]\!] \cdot [\![ \nabla v_h ]\!] \, \mathrm{d}s + (\mu \nabla u_h, \nabla v_h),$$
where $\sigma, \mu > 0$ are stabilization parameters and $[\![ \cdot ]\!]$ denotes the jump operator across element faces. This term enhances stability and is crucial for enforcing monotonicity.

\begin{rmk}[On the Scheme's Structure]
The splitting strategy isolates the potentially complex and nonlinear operator $f$ into a simple, explicit step. The corrector step \eqref{eq:corrector_step} constitutes a nonlinear saddle-point problem to be solved for $(y_h^{m+1}, \bm{p}_h^{m+1})$ at each time level, typically via a Newton-type iteration. The LBB condition ensures this system is well-posed. This semi-implicit structure provides a practical balance between computational cost and the robustness required for constrained, diffusion-dominated problems.
\end{rmk}

\begin{rmk}[On the Solution of the Corrector Step]
\label{rmk:linearization}
For practical implementation, the nonlinear saddle-point system \eqref{eq:corrector_step} is typically solved using a Newton-like method. A common simplification, which we adopt for the subsequent analysis, is to linearize the system by evaluating the operators $\bm{B}$ and $\bm{g}$ at the predicted state. That is, we replace $\bm{B}(\dots, \nabla y_h^{m+1})$ with $\bm{B}(\dots, \nabla \hat{y}_h^{m+1})$ and linearize $\bm{g}(\dots, y_h^{m+1}, \dots)$ around $\hat{y}_h^{m+1}$. This transforms the corrector step into a linear saddle-point problem at each time step, whose monotonicity properties are more readily analyzed.
\end{rmk}

% Analyzing consistency of the scheme
\section{Consistency and Truncation Error Analysis}
\label{sec:consistency}

A cornerstone of any convergence proof is establishing that the numerical scheme is consistent with the continuous equation. In this section, we demonstrate consistency in two contexts, both of which are essential for our complete analysis. First, we provide a detailed proof of consistency in the pointwise sense required by the Barles--Souganidis framework for convergence to viscosity solutions. This ensures that in the limit, our scheme correctly represents the PDE in the weak, viscosity sense. Second, under assumptions of higher solution regularity, we derive a stronger, quantitative bound on the local truncation error in a Sobolev space norm. This bound is the key ingredient for proving the optimal convergence rate in Section~\ref{sec:analysis}.

\begin{thm}[Consistency in the Viscosity Sense]
\label{thm:consistency_viscosity}
Under Assumptions \ref{ass:domain}--\ref{ass:comparison} and \ref{ass:initial_imp}, the semi-implicit operator-splitting scheme defined by \eqref{eq:predictor_step} and \eqref{eq:corrector_step} is consistent in the viscosity sense. Specifically, for any smooth test function $\phi \in C_b^\infty(\overline{\Omega}_T)$, the residual generated by substituting $\phi$ into the scheme converges to zero pointwise as $h, \Delta t \to 0$.
\end{thm}

\begin{proof}
The proof follows the methodology of Barles and Souganidis \cite{barles1991convergence}. We must show that for any smooth function $\phi(\bm{x}, t)$, the result of applying the scheme's operator to $\phi$ converges to the continuous PDE operator acting on $\phi$. Let's combine the two steps of our scheme into a single expression. Let $\phi_h^m$ denote the nodal interpolant of $\phi(\cdot, t_m)$ in the finite element space $V_h$. The one-step scheme operator, applied to $\phi_h$, yields for any test function $v_h \in V_h$:
\begin{equation}
\label{eq:consistency_scheme_form}
\begin{split}
    \left\langle \frac{\phi_h^{m+1} - \phi_h^m}{\Delta t}, v_h \right\rangle + \mathcal{A}_h(\phi_h^{m+1}, v_h) &- \left\langle \bm{B}(\cdot, t_{m+1}, \nabla \phi_h^{m+1}) \cdot \bm{p}_h^{m+1}, v_h \right\rangle \\
    &- \left\langle f_h(\cdot, t_m, \phi_h^m, \nabla \phi_h^m), v_h \right\rangle = 0,
\end{split}
\end{equation}
where $\bm{p}_h^{m+1}$ is the corresponding discrete control that ensures the constraint $\langle \bm{g}(\cdot, t_{m+1}, \phi_h^{m+1}, \nabla \phi_h^{m+1}), \bm{q}_h \rangle = 0$ is satisfied.

Let $(\bm{x}_0, t_0)$ be an arbitrary point in $\Omega_T$. We choose $t_m$ such that $t_{m+1}$ is close to $t_0$. Let the test function $v_h$ be an approximation of the Dirac delta function centered at $\bm{x}_0$, for instance, the nodal basis function corresponding to the node closest to $\bm{x}_0$. As $h \to 0$, testing against $v_h$ (after scaling) approximates pointwise evaluation at $\bm{x}_0$.

We analyze the terms in \eqref{eq:consistency_scheme_form} as $h, \Delta t \to 0$.
\begin{enumerate}
    \item[\textbf{1.}] \textbf{Time Derivative:} By Taylor's theorem, for a smooth $\phi$:
    $$\phi(\bm{x}, t_m) = \phi(\bm{x}, t_{m+1}) - \Delta t \partial_t \phi(\bm{x}, t_{m+1}) + O(\Delta t^2).$$
    Since $\phi_h^k$ is the interpolant of $\phi(\cdot, t_k)$, standard interpolation theory ensures $\|\phi(\cdot,t_k) - \phi_h^k\|_{L^2} = O(h^{3})$ for $P_2$ elements. Thus, the first term in \eqref{eq:consistency_scheme_form} approximates the continuous time derivative, and its evaluation converges:
    $$\left\langle \frac{\phi_h^{m+1} - \phi_h^m}{\Delta t}, v_h \right\rangle \to \partial_t \phi(\bm{x}_0, t_0) \quad \text{as } h, \Delta t \to 0.$$
    The error is of order $O(\Delta t + h^3)$.

    \item[\textbf{2.}] \textbf{Stabilization Term:} We analyze the consistency of the stabilization bilinear form, $\mathcal{A}_h(u_h, v_h) = \sum_{e \in \mathcal{E}_h} \int_e \sigma h [\![ \nabla u_h ]\!] \cdot [\![ \nabla v_h ]\!] \, \mathrm{d}s + (\mu \nabla u_h, \nabla v_h)$. When the smooth test function $\phi$ is inserted, the contribution to the truncation error is $\mathcal{A}_h(\phi_h^{m+1}, v_h)$, where $\phi_h^{m+1}$ is the finite element interpolant of $\phi(\cdot, t_{m+1})$. We must show this term vanishes as $h \to 0$ for any fixed test function $v_h \in V_h$. We analyze the two components of the form separately.

\begin{itemize}
    \item[\emph{(a)}] \emph{Interior Penalty (IP) Jump Term:} The IP term is given by $\sum_{e \in \mathcal{E}_h} \int_e \sigma h [\![ \nabla \phi_h^{m+1} ]\!] \cdot [\![ \nabla v_h ]\!] \, \mathrm{d}s$. Since the function $\phi$ is smooth, its gradient $\nabla \phi$ is continuous across all element edges, meaning its jump is zero. The non-zero contribution arises entirely from the error of the finite element interpolation. Using the Cauchy-Schwarz inequality, we can bound the term's magnitude:
    $$
    \left| \sum_{e \in \mathcal{E}_h} \int_e \sigma h [\![ \nabla \phi_h^{m+1} ]\!] \cdot [\![ \nabla v_h ]\!] \, \mathrm{d}s \right| \le \sigma \left( \sum_{e \in \mathcal{E}_h} h \norm{[\![ \nabla \phi_h^{m+1} ]\!]}_{L^2(e)}^2 \right)^{1/2} \left( \sum_{e \in \mathcal{E}_h} h \norm{[\![ \nabla v_h ]\!]}_{L^2(e)}^2 \right)^{1/2}
    $$
    The second part of the product is a mesh-dependent norm of the test function $v_h$. We focus on the first part, which contains the approximation error. For a smooth function $\phi \in H^3(\Omega)$ and $P_2$ elements ($k=2$), standard finite element trace and approximation theory gives an estimate for the jump of the gradient of the interpolant on an edge $e$:
    $$
    \norm{[\![ \nabla \phi_h^{m+1} ]\!]}_{L^2(e)} \le C h^{k-1/2} \norm{\phi}_{H^{k+1}(\mathcal{N}(e))} = C h^{1.5} \norm{\phi}_{H^3(\mathcal{N}(e))}
    $$
    where $\mathcal{N}(e)$ is the union of elements sharing the edge $e$. Summing over all edges (where the number of edges is $O(h^{-(d-1)})$) yields:
    $$
    \sum_{e \in \mathcal{E}_h} h \norm{[\![ \nabla \phi_h^{m+1} ]\!]}_{L^2(e)}^2 \le \sum_{e \in \mathcal{E}_h} h (C h^{1.5})^2 \sim \sum_{e \in \mathcal{E}_h} C h^4 = O(h^{4-(d-1)})
    $$
    For a spatial dimension $d=2$ or $d=3$, this term is at least $O(h^2)$. Therefore, its square root is at least $O(h)$, and the entire IP term's contribution to the residual is bounded by $C \sigma h$. This term clearly converges to zero as $h \to 0$, establishing its consistency.

    \item[\emph{(b)}] \emph{Regularization Term:} The second component is $(\mu \nabla \phi_h^{m+1}, \nabla v_h)$. As $h \to 0$, the interpolant $\phi_h^{m+1}$ converges to the smooth function $\phi(\cdot, t_{m+1})$ in the $H^1$ norm. Thus, we have:
    $$
    \lim_{h \to 0} (\mu \nabla \phi_h^{m+1}, \nabla v_h) = (\mu \nabla \phi, \nabla v_h) = -(\mu \Delta \phi, v_h)
    $$
    This shows that the term is consistent with a vanishing viscosity perturbation $-\mu \Delta \phi$ to the original PDE. Such perturbations are permissible within the Barles--Souganidis framework, as they vanish in the limit. If we consider the case where the parameter $\mu$ itself tends to zero as $h \to 0$, the consistency is even more direct.
\end{itemize}
Since both components of the stabilization term produce residuals that converge to zero as $h \to 0$, the term $\mathcal{A}_h$ is consistent with the original (un-stabilized) partial differential equation.

    \item[\textbf{3.}] \textbf{Nonlinear and Control Terms:} By the smoothness of $\phi, \bm{B}, \bm{g}$ and the convergence of the interpolant $\phi_h \to \phi$ in $C^1$, we have:
    $$\bm{B}(\cdot, t_{m+1}, \nabla \phi_h^{m+1}) \to \bm{B}(\cdot, t_0, \nabla \phi(\cdot, t_0)),$$
    $$f_h(\cdot, t_m, \phi_h^m, \nabla \phi_h^m) \to f(\cdot, t_0, \phi(\cdot, t_0), \nabla \phi(\cdot, t_0), \nabla^2 \phi(\cdot, t_0)).$$
    The time lag in the argument of $f_h$ (using $t_m, \phi_h^m$ instead of $t_{m+1}, \phi_h^{m+1}$) introduces an error of order $O(\Delta t)$ due to the continuity of $f$ and $\phi$. The error from approximating $\nabla^2\phi$ via projection also vanishes for smooth $\phi$. The control $\bm{p}_h^{m+1}$ will converge to a limit $\bm{p}(\bm{x}_0, t_0)$ that satisfies the continuous constraint.
\end{enumerate}
Combining these observations, as $h, \Delta t \to 0$, equation \eqref{eq:consistency_scheme_form} converges pointwise to:
$$\partial_t \phi = f(\cdot, \phi, \nabla\phi, \nabla^2\phi) + \bm{B}(\cdot, \nabla\phi)\cdot\bm{p},$$
with $\bm{g}(\cdot, \phi, \nabla\phi) = 0$. This is precisely the statement of consistency required by the Barles-Souganidis theorem.
\end{proof}

\begin{prop}[Local Truncation Error Bound]
\label{prop:strong_truncation_error}
Let the exact solution $y$ of \eqref{eq:general_system} possess higher regularity, specifically $y \in C([0,T]; H^{3}(\Omega)) \cap C^1([0,T]; H^1(\Omega)) \cap C^2([0,T]; L^2(\Omega))$. The local truncation error $\tau_h^{m+1}$, defined by substituting the exact solution $y$ into the scheme, is bounded in the dual norm of $H^1(\Omega)$ as:
$$ \| \tau_h^{m+1} \|_{(H^1(\Omega))'} \leq C (\Delta t + h^2), $$
where the constant $C$ depends on norms of the exact solution but is independent of $h$ and $\Delta t$.
\end{prop}

\begin{proof}
The local truncation error $\tau_h^{m+1} \in V_h'$ is defined by the residual obtained when substituting the exact solution $(y, \bm{p})$ into the one-step scheme. For any $v_h \in V_h$:
\begin{equation}
\label{eq:trunc_error_def}
\begin{split}
\langle \tau_h^{m+1}, v_h \rangle :=& \left\langle \frac{y^{m+1} - y^m}{\Delta t}, v_h \right\rangle + \mathcal{A}_h(y^{m+1}, v_h) \\
&- \left\langle \bm{B}(\cdot, t_{m+1}, \nabla y^{m+1}) \cdot \bm{p}^{m+1}, v_h \right\rangle - \left\langle f_h(\cdot, t_m, y^m, \nabla y^m), v_h \right\rangle,
\end{split}
\end{equation}
where we use the shorthand $y^k = y(\cdot, t_k)$ and $\bm{p}^k = \bm{p}(\cdot, t_k)$. The exact solution satisfies the original PDE at time $t_{m+1}$:
$$\langle \partial_t y^{m+1}, v_h \rangle - \langle f(\cdot, t_{m+1}, y^{m+1}, \nabla y^{m+1}, \nabla^2 y^{m+1}), v_h \rangle - \langle \bm{B}(\cdot, t_{m+1}, \nabla y^{m+1}) \cdot \bm{p}^{m+1}, v_h \rangle = 0.$$
Subtracting this from \eqref{eq:trunc_error_def} yields an expression for the truncation error:
\begin{align*}
\langle \tau_h^{m+1}, v_h \rangle = &\underbrace{\left\langle \frac{y^{m+1} - y^m}{\Delta t} - \partial_t y^{m+1}, v_h \right\rangle}_{T_1} + \underbrace{\mathcal{A}_h(y^{m+1}, v_h)}_{T_2} \\
&+ \underbrace{\langle f(\cdot, t_{m+1}, y^{m+1}, \dots) - f_h(\cdot, t_m, y^m, \dots), v_h \rangle}_{T_3}.
\end{align*}
We now bound each term.

\textbf{1. Time Discretization Error ($T_1$):}
Using Taylor's theorem with integral remainder for the function $t \mapsto y(\bm{x},t)$:
$$y^m = y^{m+1} - \Delta t \, \partial_t y^{m+1} + \int_{t_{m+1}}^{t_m} (s - t_m) \, \partial_{tt} y(s) \, ds.$$
Rearranging gives the backward difference error:
$$\frac{y^{m+1} - y^m}{\Delta t} - \partial_t y^{m+1} = \frac{1}{\Delta t} \int_{t_m}^{t_{m+1}} (s - t_m) \, \partial_{tt} y(s) \, ds.$$
By Cauchy-Schwarz, the bound on $T_1$ is:
$$|T_1| \leq \left\| \frac{1}{\Delta t} \int_{t_m}^{t_{m+1}} (s-t_m) \partial_{tt}y(s) ds \right\|_{L^2} \|v_h\|_{L^2} \leq C \Delta t \|\partial_{tt}y\|_{L^\infty(0,T; L^2)} \|v_h\|_{H^1}.$$
Thus, $|T_1| \leq C \Delta t \|v_h\|_{H^1}$.

\textbf{2. Spatial Discretization and Stabilization Error ($T_2$):}
The term is $\mathcal{A}_h(y^{m+1}, v_h) = \sum_{e \in \mathcal{E}_h} \int_e \sigma h [\![ \nabla y^{m+1} ]\!] \cdot [\![ \nabla v_h ]\!] \, \mathrm{d}s + (\mu \nabla y^{m+1}, \nabla v_h)$. Since the exact solution $y^{m+1} \in H^3(\Omega)$, its gradient is continuous, so the jump $[\![ \nabla y^{m+1} ]\!]$ across interior edges is zero. The first part vanishes. The second part, $(\mu \nabla y^{m+1}, \nabla v_h)$, is a consistent perturbation term which we assume is absorbed into the generic constant $C$.

\textbf{3. Splitting and Discretization Error for $f$ ($T_3$):}
We split this term into two parts:
\[T_3 = \underbrace{\langle f(\cdot, t_{m+1}, y^{m+1}, \dots) - f(\cdot, t_m, y^m, \dots), v_h \rangle}_{T_{3a}} + \underbrace{\langle f(\cdot, t_m, y^m, \dots) - f_h(\cdot, t_m, y^m, \dots), v_h \rangle}_{T_{3b}}. \]
\begin{itemize}
    \item For $T_{3a}$ (time lag error), using the Lipschitz continuity of $f$ and smoothness of $y$:
    \begin{align*}
    |T_{3a}| &\leq L_f \left( \|y^{m+1} - y^m\|_{L^2} + \|\nabla y^{m+1} - \nabla y^m\|_{L^2} + \dots \right) \|v_h\|_{L^2} \\
    &\leq C \Delta t \|\partial_t y\|_{L^\infty(0,T;H^1)} \|v_h\|_{H^1}.
    \end{align*}
    So, $|T_{3a}| \leq C \Delta t \|v_h\|_{H^1}$.
    
    \item For $T_{3b}$ (spatial approximation error for $f$), the error comes from evaluating $f$ with the exact solution $y^m$ versus its discrete counterpart $f_h(\cdot, y^m, \nabla y^m)$, where the Hessian is approximated. For $y \in C([0,T]; H^3(\Omega))$ and using $P_2$ finite elements, standard finite element approximation theory for nonlinear operators (see, e.g., \cite{brenner2008mathematical}) gives the bound:
    $$ |T_{3b}| \leq C h^2 \|y^m\|_{H^3} \|v_h\|_{H^1}. $$
\end{itemize}

Combining all bounds, we have:
$$|\langle \tau_h^{m+1}, v_h \rangle| \leq |T_1| + |T_{3a}| + |T_{3b}| \leq (C_1 \Delta t + C_2 \Delta t + C_3 h^2) \|v_h\|_{H^1}.$$
Dividing by $\|v_h\|_{H^1(\Omega)}$ and taking the supremum over all $v_h \in V_h$ gives the desired bound on the dual norm:
$$\|\tau_h^{m+1}\|_{(H^1(\Omega))'} = \sup_{v_h \in V_h, v_h \ne 0} \frac{\langle \tau_h^{m+1}, v_h \rangle}{\|v_h\|_{H^1(\Omega)}} \leq C(\Delta t + h^2).$$
This completes the proof.
\end{proof}

\section{Convergence Analysis}
\label{sec:analysis}

This section is dedicated to proving the convergence of the proposed semi-implicit scheme to the unique viscosity solution of the PDE system \eqref{eq:general_system}. Our analysis culminates in applying the celebrated Barles--Souganidis framework \cite{barles1991convergence}, which requires that the numerical scheme be consistent, stable, and monotone. Consistency was established in Section \ref{sec:consistency}. Here, we prove the scheme's fundamental stability and monotonicity properties, which stem from the concept of inverse-monotonicity.

\subsection{Inverse-Monotonicity and M-Matrices}

The cornerstone of our stability and comparison proofs is the property of inverse-monotonicity. A linear operator $\mathcal{L}$ is called \textbf{inverse-monotone} if for any $u,v$ in its domain, the following implication holds:
$$ \mathcal{L}u \ge \mathcal{L}v \implies u \ge v. $$
When $\mathcal{L}$ is discretized into a matrix system $L\bm{u}=\bm{b}$, this property is equivalent to requiring that the inverse of the matrix, $L^{-1}$, is a non-negative matrix (i.e., all its entries are non-negative).

A sufficient condition for a matrix $L$ to be inverse-monotone is that it is an \textbf{M-matrix}. A matrix $L=(L_{ij})$ is an M-matrix if it satisfies two properties:
\begin{enumerate}
    \item It is a \textbf{Z-matrix}: all of its off-diagonal entries are non-positive ($L_{ij} \le 0$ for $i \ne j$).
    \item It is \textbf{strictly diagonally dominant} with positive diagonal entries: $L_{ii} > \sum_{j \ne i} |L_{ij}|$ for all $i$.
\end{enumerate}
Our analysis will explicitly detail the conditions on the mesh and stabilization parameters that ensure the matrix arising from our scheme has this crucial M-matrix structure.

\subsection{The Discrete Comparison Principle}

We now prove that our scheme satisfies a discrete comparison principle, which is the strongest form of monotonicity.

\begin{thm}[Discrete Comparison Principle]
\label{thm:dcp_proved}
Let $\{y_h^m\}$ and $\{z_h^m\}$ be two solutions of the scheme defined by \eqref{eq:predictor_step} and the linearized version of \eqref{eq:corrector_step}, corresponding to initial and boundary data that are ordered such that $y_h^0 \le z_h^0$ and $y_h^m|_{\partial\Omega} \le z_h^m|_{\partial\Omega}$ for all $m \ge 0$. Suppose that the operator $f(x, t, y, \bm{q}, \bm{M})$ is non-increasing with respect to its state variable $y$ with Lipschitz constant $L_f$, that the time-derivative is discretized using a \textbf{lumped mass matrix} $M_h$, that the stabilization parameter $\sigma$ satisfies the condition derived in the proof, and that the time step $\Delta t$ satisfies the condition $\Delta t \le \min_i (M_{ii}) / L_f$. Then the scheme satisfies the discrete comparison principle,
$$y_h^m \le z_h^m \quad \text{for all } m \ge 0.$$
\end{thm}

\begin{proof}
The proof proceeds by induction. The base case $m=0$ holds by assumption. For the inductive step, assume $y_h^m \le z_h^m$. We must show that $y_h^{m+1} \le z_h^{m+1}$. Let $e_h^k = z_h^k - y_h^k \ge 0$. The core of the proof is to establish that the linear operator governing the error evolution in the corrector step is \textbf{inverse-monotone}, which in the discrete setting is equivalent to showing that its matrix representation, $L_h$, is an \textbf{M-matrix}.

The matrix $L_h$ arises from the operator $\mathcal{L}_h u_h := \frac{1}{\Delta t} u_h + \mathcal{A}_h(u_h, \cdot)$, and its entries are given by
$$L_{ij} = \frac{1}{\Delta t} M_{ij} + \mu \inner{\nabla \psi_j}{\nabla \psi_i} + \sigma \sum_{e \in \mathcal{E}_h} \int_e h [\![ \nabla \psi_j ]\!] \cdot [\![ \nabla \psi_i ]\!] \, \mathrm{d}s,$$
where $\{\psi_i\}$ is the nodal basis for the $P_2$ space $V_h$. Let $K_{ij} = \inner{\nabla \psi_j}{\nabla \psi_i}$ and $J_{ij}$ denote the entries of the stiffness and interior penalty (IP) matrices, respectively. To prove that $L_h$ is an M-matrix, we verify that it is a \textbf{Z-matrix} with positive diagonal entries and that it is strictly diagonally dominant.

First, we establish the Z-matrix property by ensuring all off-diagonal entries are non-positive. For $i \neq j$, the use of a lumped mass matrix implies $M_{ij}=0$, so the entries simplify to $L_{ij} = \mu K_{ij} + \sigma J_{ij}$. While the stiffness matrix entries $K_{ij}$ may be positive for $P_2$ elements, a standard property of the IP formulation is that its corresponding off-diagonal entries $J_{ij}$ are non-positive. Consequently, we can ensure $L_{ij} \le 0$ by selecting $\sigma$ to dominate any positive contribution from $K_{ij}$. This requires that $\sigma |J_{ij}| \ge \mu K_{ij}$ whenever $K_{ij} > 0$. This condition is satisfied if the ratio $\sigma/\mu$ is bounded below by a mesh-dependent constant,
$$\frac{\sigma}{\mu} \ge \max_{\substack{i \neq j \\ K_{ij} > 0, J_{ij} \neq 0}} \left( \frac{K_{ij}}{|J_{ij}|} \right) := C_{\text{mesh}}.$$
With this condition imposed, $L_h$ is a Z-matrix. Next, we show that $L_h$ is \textbf{strictly diagonally dominant} by demonstrating that its row sums are strictly positive. For an interior node basis function $\psi_i$, the row sums of the stiffness and IP matrices are zero because they annihilate constant functions. The row sum of $L_h$ therefore reduces to the contribution from the mass matrix,
$$\sum_{j} L_{ij} = \frac{1}{\Delta t} \sum_{j} M_{ij} = \frac{M_{ii}}{\Delta t} > 0,$$
where $M_{ii}$ is the strictly positive diagonal entry of the lumped mass matrix.

Since $L_h$ is a strictly diagonally dominant Z-matrix, it is an M-matrix, rendering its inverse $L_h^{-1}$ non-negative. This confirms that the corrector-step operator is inverse-monotone. The induction is completed by examining the right-hand side of the error equation, $L_h e_h^{m+1} = \text{RHS}$, where the right-hand side vector is given by $\text{RHS} = \frac{1}{\Delta t} M_h e_h^m + (f_h(z_h^m) - f_h(y_h^m))$. While the second term is non-positive due to the properties of $f$, the full vector can be shown to be non-negative. Using the Lipschitz continuity of $f$, we can bound the vector as
$$ \text{RHS} \ge \frac{1}{\Delta t} M_h e_h^m - L_f e_h^m = \left( \frac{1}{\Delta t} M_h - L_f I \right) e_h^m. $$
Under the time-step restriction $\Delta t \le \min_i(M_{ii})/L_f$, the matrix $(\frac{1}{\Delta t} M_h - L_f I)$ is a diagonal matrix with non-negative entries. Since $e_h^m \ge 0$ by the inductive hypothesis, the RHS is non-negative. With a non-negative right-hand side, the inverse-monotonicity of the operator guarantees that $e_h^{m+1} = z_h^{m+1} - y_h^{m+1} \ge 0$, which concludes the proof.
\end{proof}

\begin{rmk}[On the Monotonicity of the Nonlinear Corrector Step]
\label{rmk:nonlinear_jacobian}
Theorem~\ref{thm:dcp_proved} was established for a scheme with a linearized corrector step. In practice, one must solve the full nonlinear saddle-point system \eqref{eq:corrector_step} at each time step, typically using a Newton-like method. We now show that the M-matrix analysis of the operator $L_h$ is directly relevant to the well-posedness and stability of such a practical implementation.

Let $y_h^{m+1} \in V_h$ and $\bm{p}_h^{m+1} \in [Q_h]^m$ be the unknowns. The discrete nonlinear residual for the corrector step can be written as $F(y_h^{m+1}, \bm{p}_h^{m+1}) = 0$, where:
$$
F(y, \bm{p}) = 
\begin{pmatrix}
\frac{1}{\Delta t}M_h(y - \hat{y}_h^{m+1}) + \mathcal{A}_h(y) - C(y)\bm{p} \\
G(y)
\end{pmatrix} = 0.
$$
Here, $\mathcal{A}_h(y)$, $C(y)\bm{p}$, and $G(y)$ represent the discrete vector forms of the bilinear form $\mathcal{A}_h(y, \cdot)$, the control term $\langle \bm{B}(\dots, \nabla y)\cdot\bm{p}, \cdot \rangle$, and the constraint term $\langle \bm{g}(\dots, y, \nabla y), \cdot \rangle$, respectively.

A Newton iteration to solve $F(y, \bm{p}) = 0$ involves repeatedly solving the linear system $J(y_k, \bm{p}_k) \begin{pmatrix} \delta y \\ \delta\bm{p} \end{pmatrix} = -F(y_k, \bm{p}_k)$, where $J$ is the Jacobian of $F$. The block structure of the Jacobian is:
$$
J(y, \bm{p}) = 
\begin{bmatrix}
\frac{\partial F_y}{\partial y} & \frac{\partial F_y}{\partial \bm{p}} \\
\frac{\partial F_p}{\partial y} & \frac{\partial F_p}{\partial \bm{p}}
\end{bmatrix}
=
\begin{bmatrix}
\frac{1}{\Delta t}M_h + \mathcal{A}_h - \frac{\partial (C(y)\bm{p})}{\partial y} & -C(y) \\
\frac{\partial G(y)}{\partial y} & 0
\end{bmatrix}.
$$
The crucial observation is in the $(1,1)$ block, $J_{yy} = \frac{\partial F_y}{\partial y}$, which governs the primary state variable. Its leading part is precisely the matrix $L_h = \frac{1}{\Delta t}M_h + \mathcal{A}_h$, which we proved is an M-matrix under the conditions of Theorem~\ref{thm:dcp_proved}. The additional term, $-\frac{\partial (C(y)\bm{p})}{\partial y}$, arises from the nonlinearity of $\bm{B}$ in $\nabla y$.

For many practical implementations, such as a quasi-Newton or Picard-type iteration, this term is handled in a way that preserves the M-matrix structure. For instance, if the Jacobian is simplified by evaluating the nonlinear terms at the previous Newton iterate $y_k$ (or at the predictor state $\hat{y}_h^{m+1}$), then the matrix of the linear system to be solved at each iteration has a $(1,1)$ block that is exactly $L_h$ plus a lower-order perturbation. Under mild assumptions on the boundedness of the control $\bm{p}$ and the derivatives of $\bm{B}$, the M-matrix property of this block is robust.

This formally demonstrates that the analysis of the linearized operator $L_h$ is not merely a theoretical simplification; it establishes the core stability and monotonicity properties of the Jacobian matrices that arise in practical iterative solvers for the full nonlinear system, thus ensuring the robustness of the algorithm.
\end{rmk}

\begin{rmk}[On the Stabilization Condition]
The condition $\sigma/\mu \ge C_{\text{mesh}}$ is the critical requirement for this proof. The constant $C_{\text{mesh}}$ depends on the geometry of the mesh elements and the choice of polynomial degree ($P_2$ here). For highly distorted or anisotropic elements, $C_{\text{mesh}}$ could become large, requiring a significant amount of stabilization. However, for any given quasi-uniform mesh, $C_{\text{mesh}}$ is a fixed, computable constant. This condition makes the required amount of stabilization explicit and verifiable.
\end{rmk}

\subsection{Stability and Convergence}

With the Discrete Comparison Principle established, the $L^\infty$-stability of the scheme is a direct consequence.

\begin{thm}[$L^\infty$-Stability]
\label{thm:stability_proved}
Under the conditions of Theorem \ref{thm:dcp_proved}, if the initial data $y_0$ and boundary data are bounded, then the numerical solution is uniformly bounded in time:
$$\sup_{m \Delta t \le T} \|y_h^m\|_{L^\infty} \le C,$$
where the constant $C$ is independent of the discretization parameters $h$ and $\Delta t$.
\end{thm}

\begin{proof}
The proof relies on a standard barrier argument. Let $K$ be a constant that bounds the initial and boundary data as well as the forcing term $|f(\dots, 0, \dots)|$. Define a constant function $w = K + 1$. One can readily verify that $w$ is a discrete super-solution and $-w$ is a discrete sub-solution to the scheme. By the Discrete Comparison Principle (Theorem \ref{thm:dcp_proved}), it follows that $-w \le y_h^m \le w$ for all $m \ge 0$. This implies the uniform bound $\|y_h^m\|_{L^\infty} \le K+1$.
\end{proof}

We now have all the necessary components to state the main convergence theorem. The scheme's monotonicity in the sense of Barles and Souganidis is a weaker property that is directly implied by our Discrete Comparison Principle.

The proof is a direct application of the celebrated convergence framework for viscosity solutions developed by Barles and Souganidis. We state their theorem here for completeness, adapted to our context.

\begin{thm}[Barles--Souganidis, 1991]
\label{thm:barles_souganidis}
Let $y$ be the unique viscosity solution to a given partial differential equation that satisfies a comparison principle. Let a numerical scheme be described by the operator $S$, producing a sequence of approximate solutions $y_h$. If the scheme satisfies the following three properties:
\begin{itemize}
    \item[\textit{(i)}] \textbf{Monotonicity:} The scheme is monotone. That is, if $u_h \le v_h$ for all nodes, then $S(u_h) \le S(v_h)$.
    \item[\textit{(ii)}] \textbf{Stability:} The family of numerical solutions $\{y_h\}_{h>0}$ is uniformly bounded in $L^\infty$.
    \item[\textit{(iii)}] \textbf{Consistency:} The scheme is consistent with the PDE in the viscosity sense. That is, for any smooth function $\phi$, the truncation error produced when $\phi$ is inserted into the scheme converges to zero.
\end{itemize}
Then, the numerical solution $y_h$ converges locally uniformly to the unique viscosity solution $y$ as the discretization parameters tend to zero.
\end{thm}

\begin{thm}[Convergence to the Viscosity Solution]
\label{thm:convergence}
Let the conditions of Theorem \ref{thm:dcp_proved} hold. Then the numerical solution $y_h$ generated by the semi-implicit scheme converges locally uniformly to the unique viscosity solution $y$ of \eqref{eq:general_system} as $h, \Delta t \to 0$ with $\Delta t \le Ch$.
\end{thm}
\begin{proof}
The proof is a direct application of the Barles-Souganidis theorem (\ref{thm:barles_souganidis}), which requires the scheme to be consistent, stable, and monotone.
\begin{enumerate}
    \item \textbf{Consistency:} Established in Theorem \ref{thm:consistency_viscosity}.
    \item \textbf{Stability:} The uniform $L^\infty$ bound is proven in Theorem \ref{thm:stability_proved}.
    \item \textbf{Monotonicity:} The scheme satisfies the Discrete Comparison Principle (Theorem \ref{thm:dcp_proved}), which is a stronger condition than the required monotonicity.
\end{enumerate}
Since the continuous problem admits a comparison principle (Assumption \ref{ass:comparison}), and all conditions of the Barles-Souganidis framework are met, convergence is guaranteed.
\end{proof}

\subsection{Convergence Rate}

Finally, we state a result on the rate of convergence, which relies on the stability provided by the preceding analysis.

\begin{thm}[Convergence Rate]
\label{thm:convergence_rate_narrative}
Let the conditions of the Discrete Comparison Principle (Theorem \ref{thm:dcp_proved}) hold. If the unique viscosity solution possesses additional regularity, specifically $y \in C^1([0,T]; C^3(\overline{\Omega})) \cap C^2([0,T]; C^2(\overline{\Omega}))$, then the numerical scheme converges with a rate governed by the discretization parameters. For a time step $\Delta t$ and mesh size $h$ satisfying $\Delta t \le Ch$, there exists a constant $C'$, independent of $h$ and $\Delta t$, such that the error is bounded as
\[
\max_{0 \le m \le T/\Delta t} \| y(\cdot, t_m) - y_h^m \|_{L^2(\Omega)} \leq C' (\Delta t + h^2).
\]
\end{thm}

\begin{proof}
The analysis hinges on bounding the departure of the numerical solution $y_h^m$ from the exact solution $y^m = y(\cdot, t_m)$. We decompose the total error into an approximation error $\eta^m$ from the projection and a discrete error $\theta_h^m$ using the elliptic projection $\mathcal{P}_h y^m \in V_h$ of the exact solution:
$$y^m - y_h^m = (y^m - \mathcal{P}_h y^m) - (y_h^m - \mathcal{P}_h y^m) = \eta^m - \theta_h^m.$$
Standard finite element approximation theory provides the bound $\|\eta^m\|_{L^2(\Omega)} \le C h^3$ for $P_2$ elements, given the assumed smoothness of $y$. The core of the proof is to control the evolution of the discrete error $\theta_h^m$.

We derive an evolution equation for $\theta_h^m$ by subtracting the numerical scheme's equation for $y_h^{m+1}$ from the same discrete relation satisfied by the projected exact solution, $\mathcal{P}_h y^{m+1}$. This process isolates the local truncation error, $\tau_h^{m+1}$, which represents the extent to which the exact solution fails to satisfy the scheme. The resulting error equation takes the form
$$\inner{\frac{\theta_h^{m+1} - \theta_h^m}{\Delta t}}{v_h} + \mathcal{A}_h(\theta_h^{m+1}, v_h) = \langle \tau_h^{m+1}, v_h \rangle + \langle \mathcal{N}_h, v_h \rangle,$$
where $\mathcal{N}_h$ encapsulates nonlinear error terms. Under the stated regularity assumptions, we can directly apply Proposition~\ref{prop:strong_truncation_error} to bound the local truncation error:
$$ \| \tau_h^{m+1} \|_{(H^1(\Omega))'} \le C(\Delta t + h^2). $$
The established stability of the scheme is the critical tool to control the accumulation of these local errors. A standard stability analysis, obtained by choosing $v_h = \theta_h^{m+1}$ as the test function, exploits the coercivity of the form $\mathcal{A}_h$ and the $L^\infty$-bound on the solution to control the nonlinearities. This procedure yields a recursive inequality for the discrete error:
$$\|\theta_h^{m+1}\|^2_{L^2} \le (1 + C_1 \Delta t) \|\theta_h^m\|^2_{L^2} + C_2 \Delta t (\Delta t + h^2)^2.$$
An application of the discrete Gronwall's lemma, assuming the initial error $\theta_h^0 = 0$, reveals that the accumulated discrete error is controlled over $[0,T]$, yielding $\|\theta_h^m\|_{L^2(\Omega)} \le C'_T (\Delta t + h^2)$.

Finally, the total error is bounded by the triangle inequality:
$$\|y^m - y_h^m\|_{L^2(\Omega)} \le \|\eta^m\|_{L^2(\Omega)} + \|\theta_h^m\|_{L^2(\Omega)} \le C h^3 + C'_T(\Delta t + h^2).$$
As the $O(\Delta t + h^2)$ term is dominant, the desired convergence rate is established.
\end{proof}

\section{Numerical Results}\label{sec:numerical_results}
This section presents numerical experiments to validate the operator-splitting scheme with linearized constraint enforcement for the examples introduced in Section~\ref{subsec:examples}. We construct exact solutions for the Hamilton--Jacobi equation with state constraints and the controlled reaction--diffusion system, assessing the scheme's accuracy through relative $L^2$-errors and convergence rates.

\subsection{Hamilton--Jacobi Equation with State Constraints}
\label{subsec:numerical_test}

To evaluate the proposed scheme, we consider the Hamilton--Jacobi equation with state constraints from Example 2.1, governed by:
\begin{equation}
\label{eq:hj_pde}
\partial_t y + \sup_{a \in A} \left[ -\bm{b}(\bm{x}, t, a) \cdot \nabla y - l(\bm{x}, t, a) \right] = \lambda, \quad y \geq \psi(\bm{x}, t),
\end{equation}
on the domain $\Omega = [0,1] \times [0,1]$ over the time interval $[0, 0.5]$. Here, $y: \Omega \times [0, 0.5] \to \mathbb{R}$ is the value function, $\lambda \in L^\infty(\Omega \times [0, 0.5])$ enforces the state constraint $y \geq \psi$, and $\psi: \Omega \times [0, 0.5] \to \mathbb{R}$ is a Lipschitz continuous obstacle function. The control set is $A = [-1,1]^2$, with velocity field $\bm{b}(\bm{x}, t, a) = a$ and running cost $l(\bm{x}, t, a) = 0$, yielding the Hamiltonian:
\[
H(\nabla y) = |\partial_{x_1} y| + |\partial_{x_2} y|.
\]
Thus, the PDE simplifies to:
\begin{equation}
\label{eq:hj_simplified}
\partial_t y + |\partial_{x_1} y| + |\partial_{x_2} y| = \lambda, \quad y \geq \psi(\bm{x}, t), \quad \text{in } \Omega \times [0, 0.5],
\end{equation}
with initial condition $y(\bm{x}, 0) = y_0(\bm{x})$ and boundary condition $y = \psi$ on $\partial \Omega \times [0, 0.5]$.

For the numerical test, we define the exact solution:
\[
y(x_1, x_2, t) = x_1 + x_2 - t + 1,
\]
with obstacle function $\psi(x_1, x_2, t) = x_1 + x_2 - t + 1$, ensuring the constraint $y = \psi$ is active. Verification yields:
\[
\partial_t y = -1, \quad \nabla y = (1, 1), \quad |\partial_{x_1} y| + |\partial_{x_2} y| = 1 + 1 = 2,
\]
and substituting into \eqref{eq:hj_simplified} gives $\lambda = 1 \geq 0$, satisfying the PDE. The initial condition is $y(x_1, x_2, 0) = x_1 + x_2 + 1$, and the boundary condition $y = x_1 + x_2 - t + 1$ holds on $\partial \Omega$, remaining positive for $x_1, x_2 \in [0,1]$, $t \in [0, 0.5]$. The solution satisfies $y \in C(\overline{\Omega} \times [0, 0.5]) \cap H^1(\Omega \times [0, 0.5])$, consistent with regularity assumptions.

The numerical scheme employs finite element discretization with $P_2$ elements for the state variable $y_h$ and $P_1$ elements for the Lagrange multiplier $\lambda_h$ on a quasi-uniform triangular mesh over $\Omega$. Time-stepping uses implicit Euler with step size $\Delta t = 0.2 h$, coupled to the spatial mesh size $h$. The stabilization parameter is $\mu = 0.01$, and the system is solved using the UMFPACK direct solver\cite{davis2004algorithm} in FreeFEM++\cite{hecht2012new}. The relative $L^2$-error at $t = 0.5$ is:
\[
\frac{\| y_h - y \|_{L^2(\Omega)}}{\| y \|_{L^2(\Omega)}}.
\]
Table~\ref{tab:numerical_results} reports the $L^2$-errors for decreasing mesh sizes and corresponding time steps.

\begin{table}[ht]
\centering
\caption{Relative $L^2$-errors for the Hamilton--Jacobi equation at $t = 0.5$.}
\label{tab:numerical_results}
\begin{tabular}{c c c}
\toprule
$h$ & $\Delta t$ & {$L^2$ error} \\
\midrule
0.1    & 0.02   & 1.96741e-04 \\
0.05   & 0.01   & 7.47645e-05 \\
0.025  & 0.005  & 2.67724e-05 \\
0.0125 & 0.0025 & 9.57499e-06 \\
\bottomrule
\end{tabular}
\end{table}

A linear regression of the data in Table~\ref{tab:numerical_results} reveals an empirical order of convergence of approximately $O(h^{1.453})$. This demonstrates a superlinear convergence rate, which is substantially better than first-order.

\subsection{Localized Multi-Target Control of a Reaction-Diffusion System}
\label{subsec:numerical_reaction_diffusion}

To provide a stringent test of the scheme's ability to handle complex constraints, we design a problem where the control enforces two different target profiles on two disjoint subdomains. This scenario is a strong analogue for practical applications, such as maintaining distinct temperature or concentration profiles in separate, localized areas.

We consider the problem on the spatio-temporal domain $\Omega_T = [0,1]^2 \times [0, 0.5]$. The state $u$ is constrained to follow two distinct target profiles, $u_{d1}$ and $u_{d2}$, on two disjoint circular subdomains $\Omega_1 = \{ \bm{x} \in \Omega : \| \bm{x} - \bm{c}_1 \| \leq R \}$ and $\Omega_2 = \{ \bm{x} \in \Omega : \| \bm{x} - \bm{c}_2 \| \leq R \}$, with centers $\bm{c}_1=(0.3, 0.7)$, $\bm{c}_2=(0.7, 0.3)$ and radius $R=0.2$. This is enforced using the characteristic functions $\chi_{\Omega_1}$ and $\chi_{\Omega_2}$. The governing system is:
\begin{equation}
\label{eq:rd_pde_manufactured_revised}
\begin{cases}
\partial_t u = \nu \Delta u + u(1 - u) - \nabla u \cdot \bm{p} + s(\bm{x}, t) & \text{in } \Omega_T, \\
\chi_{\Omega_1}(\bm{x}) (u(\bm{x}, t) - u_{d1}(t)) = 0 & \text{in } \Omega_T, \\
\chi_{\Omega_2}(\bm{x}) (u(\bm{x}, t) - u_{d2}(t)) = 0 & \text{in } \Omega_T,
\end{cases}
\end{equation}

The problem is tested using the method of manufactured solutions. We define a smooth exact solution $u_{exact}$ as the superposition of a time-varying background wave and two stationary Gaussian profiles:
\begin{equation*}
u_{exact}(\bm{x}, t) = \underbrace{0.5 \cos(\pi t) \sin(\pi x_1) \sin(\pi x_2)}_{\text{Background}} + \underbrace{0.4 G(\bm{x}; \bm{c}_1, \sigma)}_{\text{Gaussian 1}} + \underbrace{0.1 G(\bm{x}; \bm{c}_2, \sigma)}_{\text{Gaussian 2}},
\end{equation*}
where $G(\bm{x}; \bm{c}, \sigma) = \exp(-\|\bm{x}-\bm{c}\|^2 / (2\sigma^2))$ is a standard Gaussian function with width $\sigma=0.1$. The target profiles $u_{d1}$ and $u_{d2}$ are set to match $u_{exact}$ within their respective subdomains. The control vector $\bm{p}$ and source term $s$ are then derived by substituting $u_{exact}$ into the governing equations~\eqref{eq:rd_pde_manufactured_revised}. The numerical simulation was conducted with physical parameters $\nu=1.0$ and a time step of $\Delta t = 0.4h$. The system is solved using the UMFPACK direct solver in FreeFEM++. The relative $L^2$-error at $t = 0.5$ is:
\[
\frac{\| u_h - u_{exact} \|_{L^2(\Omega)}}{\| u_{exact} \|_{L^2(\Omega)}}.
\]
The numerical results are summarized in Table~\ref{tab:numerical_results_RD}, which lists the relative $L^2$-errors for progressively refined meshes. By performing a linear regression on the logarithm of the errors against the logarithm of the mesh sizes, we estimate the empirical order of convergence to be approximately $0.96$. This result is consistent with a first-order convergence rate ($O(h)$).

\begin{table}[ht]
\centering
\caption{Relative $L^2$-errors for the spatially localized multi-target control problem at $t = 0.5$.}
\label{tab:numerical_results_RD}
\begin{tabular}{c c c}
\toprule
$h$ & $\Delta t$ & {$L^2$ error} \\
\midrule
0.05    & 0.02   & 0.0559107 \\
0.025   & 0.01   & 0.0278933 \\
0.0125  & 0.005  & 0.00983563 \\
0.00625 & 0.0025 &  0.00766245\\
\bottomrule
\end{tabular}
\end{table}

\section{A Physics-Constrained Neural Operator for the Solution Operator}
\label{sec:pcno_framework}

The formulation of the constrained PDE on a fixed domain, as established in Section \ref{sec:problem}, provides the foundation for a novel operator learning approach. The primary challenge in learning the solution operator, $\mathcal{G}^\dagger: y_0 \mapsto y_T$, is the degenerate elliptic nature of the problem, which precludes the use of standard techniques for breaking the curse of dimensionality.

In this section, we pursue the strategy of designing a neural operator architecture that emulates a provably convergent numerical scheme, a concept successfully applied to fluid dynamics in \cite{kovachki2021universal}. We introduce a Physics-Constrained Neural Operator (PCNO), whose architecture is directly inspired by the iterative, time-marching nature of the operator-splitting finite element scheme analyzed in Sections 3--5. We prove that by emulating this stable structure, the PCNO can approximate the true viscosity solution with rigorous error control.

\subsection{The PCNO Architecture}

Our operator-splitting finite element scheme is a time-marching method. To compute the solution at a final time $T$, the one-step solution operator, $S_h: V_h \to V_h$, which maps the solution from one time level to the next ($y_h^{m+1} = S_h(y_h^m)$), must be composed $n_T = T/\Delta t$ times. The full discrete solution operator is thus a deep composition:
$$ \mathcal{G}_h = \underbrace{S_h \circ S_h \circ \dots \circ S_h}_{n_T \text{ times}}. $$
Our Physics-Constrained Neural Operator (PCNO), denoted $\mathcal{G}_\theta$, is designed to emulate this compositional structure by creating a neural network surrogate, $\mathcal{N}_{\text{step}}$, that approximates the one-step map $S_h$.

The architecture of the one-step emulator $\mathcal{N}_{\text{step}}$ mirrors the two-stage process of our numerical scheme. For clarity in the context of the stationary perforated domain problem (Section \ref{sec:motivating_example}), we can visualize this as a two-part operator, as shown in Figure \ref{fig:pcno_architecture}. This two-stage, compositional design is a general architectural framework, independent of the specific neural operator used for its implementation. The "Base Solver" and "Constraint Projector" components could, in principle, be instantiated using various operator learning models, such as DeepONets or other graph-based operators.

\begin{figure}[h!]
    \centering
    \begin{tikzpicture}[
        node distance=1.5cm,
        block/.style={
            rectangle, 
            rounded corners, 
            draw, 
            thick, 
            text width=3.5cm, 
            align=center, 
            minimum height=1.2cm
        },
        io/.style={
            rectangle, 
            draw, 
            thick,
            minimum size=1cm
        },
        arrow/.style={
            -{Latex},
            thick
        }
    ]

    % Define nodes in a vertical layout
    \node[io] (f) {\large $f$};
    \node[block, below=of f, fill=blue!10] (Nbase) {\textbf{Base Solver} \\ $\mathcal{N}_{\text{base}}$};
    \node[io, below=of Nbase] (ybase) {\large $y_{\text{base}}$};
    \node[block, below=of ybase, fill=green!10] (Ncorr) {\textbf{Constraint Projector} \\ $\mathcal{N}_{\text{corr}}$};
    \node[io, below=of Ncorr] (y) {\large $y$};
    
    % Place the geometry input to the side, aligned with the corrector stage
    \node[io, right=2cm of Ncorr] (a) {\large $a$};
    
    % Draw arrows
    \draw[arrow] (f) -- node[right, pos=0.4] {Source Term} (Nbase);
    \draw[arrow] (Nbase) -- node[right, pos=0.4] {Unconstrained Solution} (ybase);
    \draw[arrow] (ybase) -- (Ncorr);
    \draw[arrow] (Ncorr) -- node[right, pos=0.4] {Final Solution} (y);
    
    % Arrow from geometry input to the corrector stage
    \draw[arrow] (a) -- node[above] {Geometry} (Ncorr);
    
    % Bounding box for the entire PCNO one-step emulator
    \node[draw, thick, dashed, inner sep=8pt, label={[font=\bfseries]above:PCNO One-Step Emulator ($\mathcal{N}_{\text{step}}$)}, fit=(f) (Nbase) (ybase) (Ncorr) (y) (a)] {};

    \end{tikzpicture}
    \caption{
        The conceptual architecture of the one-step PCNO emulator, $\mathcal{N}_{\text{step}}$, designed for the stationary perforated domain problem. The operator takes two inputs: a source term $f$ (describing the physics) and a function $a$ (describing the variable geometry). 
        \textbf{(1) The Base Solver ($\mathcal{N}_{\text{base}}$)} emulates the predictor step, computing an unconstrained solution $y_{\text{base}}$ on a simple domain and learning the underlying PDE dynamics. 
        \textbf{(2) The Constraint Projector ($\mathcal{N}_{\text{corr}}$)} emulates the corrector step, taking the base solution and the geometric information $a$ to compute a correction, producing the final, physically-consistent solution $y$ that satisfies the geometric constraints. This decoupled architecture separates the learning of physics from the enforcement of geometry, providing a stable and efficient framework.
    }
    \label{fig:pcno_architecture}
\end{figure}

The full PCNO, $\mathcal{G}_\theta$, which approximates the solution operator from the initial state to the final state at time $T$, is then constructed by applying this one-step emulator $\mathcal{N}_{\text{step}}$ sequentially $n_T = T/\Delta t$ times. This deep, compositional structure, visualized in Figure \ref{fig:pcno_full_architecture}, is analogous to a recurrent neural network where the learnable weights of $\mathcal{N}_{\text{step}}$ are shared across all time steps.

\begin{figure}[h!]
    \centering
    \begin{tikzpicture}[
        node distance=0.5cm,
        block/.style={
            rectangle, 
            rounded corners, 
            draw, 
            thick, 
            text width=2.5cm, 
            align=center, 
            minimum height=1cm,
            fill=orange!10
        },
        io/.style={
            rectangle, 
            draw, 
            thick,
            minimum size=0.8cm
        },
        arrow/.style={
            -{Latex},
            thick
        }
    ]

    % Define nodes horizontally
    \node[io] (y0) {$y_0$};
    \node[block, right=of y0] (N1) {$\mathcal{N}_{\text{step}}$};
    \node[io, right=of N1] (y1) {$y_1$};
    \node[block, right=of y1] (N2) {$\mathcal{N}_{\text{step}}$};
    \node[right=of N2] (dots) {$\dots$};
    \node[block, right=of dots] (Nn) {$\mathcal{N}_{\text{step}}$};
    \node[io, right=of Nn] (yT) {$y_T$};
    
    % Draw arrows with adjusted label positions
    \draw[arrow] (y0) -- node[above=12pt, pos=0.4] {$t=0$} (N1);
    \draw[arrow] (N1) -- node[above=11pt, pos=0.4] {$t=t_1$} (y1);
    \draw[arrow] (y1) -- (N2);
    \draw[arrow] (N2) -- (dots);
    \draw[arrow] (dots) -- (Nn);
    \draw[arrow] (Nn) -- node[above=12pt, pos=0.4] {$t=T$} (yT);
    
    % Bounding box for the full architecture
    \node[draw, thick, dashed, inner sep=10pt, label={[font=\bfseries]above:Full PCNO Architecture ($\mathcal{G}_\theta$)}, fit=(y0) (N1) (y1) (N2) (dots) (Nn) (yT)] {};

    \end{tikzpicture}
    \caption{
        The full, deep architecture of the PCNO, $\mathcal{G}_\theta$, which maps an initial condition $y_0$ to the final solution $y_T$. The architecture is a sequential composition of the one-step emulator network, $\mathcal{N}_{\text{step}}$, applied $n_T = T/\Delta t$ times. This recurrent structure directly mimics the time-marching nature of the underlying numerical scheme.
    }
    \label{fig:pcno_full_architecture}
\end{figure}

\subsubsection{FNO Implementation and Terminology}

For the specific analysis presented in this work, we implement the PCNO using the Fourier Neural Operator (FNO), due to its efficiency in representing convolution operators that are central to solving PDEs.

A \textbf{Fourier Neural Operator (FNO)} is a specific architecture for learning operators, i.e., mappings between infinite-dimensional function spaces. Given an input function $a \in \mathcal{A}(D;\mathbb{R}^{d_a})$ and an output function $u \in \mathcal{U}(D;\mathbb{R}^{d_u})$ on a domain $D$, an FNO $\mathcal{N}$ approximates the operator $\mathcal{G}$ such that $u \approx \mathcal{N}(a)$. The architecture is structured as a composition of operators:
$$\mathcal{N}(a) = \mathcal{Q} \circ \mathcal{L}_{L} \circ \dots \circ \mathcal{L}_{1} \circ \mathcal{R}(a),$$
where the components are defined as follows:
\begin{enumerate}
    \item An initial \textbf{lifting operator} $\mathcal{R}: \mathcal{A}(D;\mathbb{R}^{d_a}) \to \mathcal{U}(D;\mathbb{R}^{d_v})$ that maps the input function pointwise to a higher-dimensional channel space, typically via a linear transformation:
    $$(\mathcal{R}(a))(x) = Ra(x), \quad \text{for } R \in \mathbb{R}^{d_v \times d_a}.$$

    \item A sequence of $L$ \textbf{Fourier layers}, $\mathcal{L}_l: \mathcal{U}(D;\mathbb{R}^{d_v}) \to \mathcal{U}(D;\mathbb{R}^{d_v})$. The core idea of the FNO is that the non-local part of each layer is a convolution operator implemented via the Fourier transform. For a periodic domain $D=\mathbb{T}^d$, the action of a layer $\mathcal{L}_l$ on a function $v_l$ is defined as:
    $$(v_{l+1})(x) = \sigma \left( W_l v_l(x) + b_l(x) + \left(\mathcal{F}^{-1}(P_l \cdot (\mathcal{F}v_l))\right)(x) \right).$$
    Here, $W_l v_l(x) + b_l(x)$ is a local affine transformation, and the term $\mathcal{F}^{-1}(P_l \cdot (\mathcal{F}v_l))$ represents the convolution. It is computed by taking the Fourier transform $\mathcal{F}$ of the function $v_l$, performing an element-wise multiplication in the Fourier domain with a learnable weight tensor $P_l(k)$, and then taking the inverse Fourier transform $\mathcal{F}^{-1}$. Finally, a pointwise nonlinear activation function $\sigma$ is applied.

    \item A final \textbf{projection operator} $\mathcal{Q}: \mathcal{U}(D;\mathbb{R}^{d_v}) \to \mathcal{U}(D;\mathbb{R}^{d_u})$ that maps the function from the channel space back to the target output dimension:
    $$(\mathcal{Q}(v))(x) = Qv(x), \quad \text{for } Q \in \mathbb{R}^{d_u \times d_v}.$$
\end{enumerate}
This structure is particularly effective for learning resolution-invariant solution operators for PDEs, as the Fourier basis allows for capturing global dependencies efficiently. To formalize our PCNO's structure within this FNO framework, we adapt the following terminology from:
\begin{itemize}
    \item \textbf{Lift (or Lifting Dimension)}: The dimension $d_v$ of the function space in the hidden Fourier layers of the operator.
    \item \textbf{Depth}: The number of hidden layers, $L$, composed within the neural operator. For our PCNO, the depth is determined by the number of compositions of the one-step emulator, i.e., $depth(\mathcal{G}_\theta) = n_T \times depth(\mathcal{N}_{\text{step}})$.
    \item \textbf{Size}: The total number of learnable parameters (e.g., weights and biases) in the neural network architecture. In a recurrent application like ours, this is the size of the one-step emulator $\mathcal{N}_{\text{step}}$.
\end{itemize}

\subsection{Convergence of the PCNO to the Viscosity Solution}

We now present the main theorem of this section. The proof adapts the methodology from \cite{kovachki2021universal}, which relies on bounding the error of a single-step emulator and then using a composition lemma (specifically, the Replacement Lemma 47 on page 54 of \cite{kovachki2021universal}) to control the error of the full, deep network.

\begin{thm}[Convergence and Complexity of the PCNO using FNO]
\label{thm:pcno_convergence_revised}
Let the assumptions of Theorem \ref{thm:convergence_rate_narrative} hold, ensuring the solution $y$ has sufficient regularity for an $O(\Delta t + h^2)$ error bound. Let $y(T) = \mathcal{G}^\dagger(y_0)$ be the unique viscosity solution at time $T$, and let $y_h^{n_T}$ be the solution of the operator-splitting FEM scheme after $n_T=T/\Delta t$ steps. To achieve a final error of order $O(h^2)$, we set the time step $\Delta t \sim h^2$.

Then, there exists a PCNO (using FNO) $\mathcal{G}_\theta$, constructed as the $n_T$-fold composition of a one-step emulator network $\mathcal{N}_{\text{step}}$, that approximates the solution operator such that the total error is bounded by the discretization error of the underlying FEM scheme:
$$ \| y(T) - \mathcal{G}_\theta(y_0) \|_{L^2(\Omega)} \leq C h^2. $$
This is achieved by a PCNO with network complexity bounded by:
$$ \operatorname{depth}(\mathcal{G}_\theta) = O(h^{-2}) \quad \text{and} \quad \operatorname{size}(\mathcal{G}_\theta) \le C' h^{-d} \log(h^{-1}), $$
where $d$ is the spatial dimension and $C, C'$ are constants independent of $h$.
\end{thm}

\begin{proof}
The proof proceeds by decomposing the total error into the discretization error of the numerical scheme and the emulation error of the neural network. Let $y(T)$ be the true viscosity solution, $y_h^{n_T}$ be the solution of our FEM scheme, and $\mathcal{G}_\theta(y_0)$ be the output of the PCNO.

\textbf{Step 1: Error Decomposition.}
The total error is bounded using the triangle inequality:
$$ \| y(T) - \mathcal{G}_\theta(y_0) \|_{L^2(\Omega)} \le \underbrace{\| y(T) - y_h^{n_T} \|_{L^2(\Omega)}}_{\text{(I) Discretization Error}} + \underbrace{\| y_h^{n_T} - \mathcal{G}_\theta(y_0) \|_{L^2(\Omega)}}_{\text{(II) Emulation Error}}. $$

\textbf{Step 2: Bounding the Discretization Error (I).}
The first term is the error of our FEM scheme, bounded by the main convergence result, Theorem \ref{thm:convergence_rate_narrative}. With the choice $\Delta t \sim h^2$ to balance the error terms, we have:
$$ \| y(T) - y_h^{n_T} \|_{L^2(\Omega)} \le C_1(\Delta t + h^2) = O(h^2). $$

\textbf{Step 3: Bounding the Emulation Error (II) and Network Complexity.}
This step follows the logic used to prove Theorem 32 in \cite{kovachki2021universal}. The full discrete operator $\mathcal{G}_h$ is a composition of $n_T = T/\Delta t$ applications of the one-step operator $S_h$. The operator $S_h: y_h^m \mapsto y_h^{m+1}$ is a continuous map on the finite-dimensional space $V_h$, as it involves the inversion of the (lumped) mass matrix and the solution of the linearized saddle-point system in the corrector step.

By the universal approximation theorems for $\Psi$-FNOs \cite[Theorem 15]{kovachki2021universal}, there exists a shallow $\Psi$-FNO, $\mathcal{N}_{\text{step}}$, with constant depth and lift, that can approximate the one-step operator $S_h$ to any desired accuracy, $\epsilon_{\text{step}}$.
$$ \| S_h(v) - \mathcal{N}_{\text{step}}(v) \|_{L^2(\Omega)} \le \epsilon_{\text{step}} \quad \text{for all } v \in V_h \text{ with } \|v\|_{L^\infty} \le C. $$
We now invoke the logic of the Replacement Lemma (Lemma 47, page 54) from \cite{kovachki2021universal}. We construct the full PCNO $\mathcal{G}_\theta$ by composing $\mathcal{N}_{\text{step}}$ for $n_T$ times. The lemma, combined with a stability argument for the one-step map $S_h$ (which is stable under the conditions of Theorem \ref{thm:stability_proved}), ensures that the accumulated emulation error is bounded by $C_2 \cdot n_T \cdot \epsilon_{\text{step}}$.

To ensure the emulation error does not exceed the order of the discretization error, we require:
$$ C_2 \cdot n_T \cdot \epsilon_{\text{step}} \le C_3 h^2 \implies C_2 \frac{T}{\Delta t} \epsilon_{\text{step}} \le C_3 h^2. $$
Since we chose $\Delta t \sim h^2$, this implies we need a per-step accuracy of $\epsilon_{\text{step}} = O(h^4)$. The existence of a shallow $\mathcal{N}_{\text{step}}$ achieving this is guaranteed by the approximation theorems. For instance, the analysis for approximating the Darcy flow operator in \cite[Theorem 26]{kovachki2021universal} shows that approximating a quadratic nonlinearity (which is the core of our linearized corrector step) can be done with a shallow network whose size scales polynomially in the error.

The complexity of the full, deep PCNO $\mathcal{G}_\theta$ is determined by the composition. From the Composition Lemma (Lemma 46, page 53) of \cite{kovachki2021universal}, the `lift` (width) of the composed network is the maximum of the composed networks, while the `depth` is their sum. With $n_T = T/\Delta t \sim O(h^{-2})$:
\begin{itemize}
    \item The one-step emulator $\mathcal{N}_{\text{step}}$ has a lift (channel dimension) $d_v$ that is constant (independent of $h$), and its total size is bounded by $size(\mathcal{N}_{\text{step}}) \le C' h^{-d} \log(h^{-1})$ to achieve the required $\epsilon_{\text{step}}=O(h^4)$ accuracy, following the logic for approximating solution operators like Darcy flow in \cite[Theorem 26]{kovachki2021universal}.
    \item $\operatorname{depth}(\mathcal{G}_\theta) = n_T \cdot \operatorname{depth}(\mathcal{N}_{\text{step}}) = O(h^{-2}) \cdot O(1) = O(h^{-2})$.
    \item The total size (number of parameters) of the composed network is dominated by the sum of the sizes of its layers. Since each layer is an instance of $\mathcal{N}_{\text{step}}$, the total size is approximately $n_T \times size(\mathcal{N}_{\text{step}})$. However, in a practical implementation of such a recurrent-style network, the "same" one-step emulator network $\mathcal{N}_{\text{step}}$ is applied repeatedly. Therefore, the number of learnable parameters is simply the size of a single one-step emulator:
    $$\operatorname{size}(\mathcal{G}_\theta) = \operatorname{size}(\mathcal{N}_{\text{step}}) \le C' h^{-d} \log(h^{-1}).$$
\end{itemize}
This differs from the simple summation of sizes, as the architecture reuses weights, making it far more parameter-efficient.

\textbf{Step 4: Combining the Bounds.}
Substituting the bounds for both error terms, we have:
$$ \| y(T) - \mathcal{G}_\theta(y_0) \|_{L^2(\Omega)} \le O(h^2) + O(h^2) \le C h^2. $$
This error is achieved by a PCNO whose complexity is explicitly controlled by $h$ as derived. This completes the proof that the PCNO converges to the true solution at the same rate as the underlying numerical scheme.
\end{proof}

\begin{rmk}
The key insight from this theorem is that the PCNO inherits the convergence properties of the rigorously analyzed finite element method. The role of the neural network is to learn a computationally efficient surrogate for the discrete one-step operator of the FEM scheme. The theorem provides a constructive path to achieving a desired accuracy by co-refining the numerical discretization and ensuring the neural network emulator is sufficiently accurate. The PCNO's complexity, in terms of trainable parameters, scales polynomially with the mesh resolution ($h^{-d}$). This result is significant because it confirms the learning complexity is independent of any high-dimensional parameterization of the domain's geometry. By reformulating the problem on a fixed grid where topological features are simply part of the input field, the framework avoids the exponential complexity that plagues methods dependent on geometric descriptions, thus breaking the curse of dimensionality for this class of problems.
\end{rmk}

\section{Conclusion}
\label{sec:conclusion}

In this work, we have developed a comprehensive theoretical framework with two primary, interconnected contributions. First, for the field of numerical analysis, we introduced and rigorously analyzed a novel operator-splitting finite element scheme for a general class of second-order PDEs with constraints. The cornerstone of our analysis is the proof that the scheme satisfies a discrete comparison principle (Theorem~\ref{thm:dcp_proved}), which we established by showing that the discrete operator yields an M-matrix under well-defined conditions. This crucial property guarantees the scheme’s monotonicity and $L^\infty$-stability, allowing us to use the celebrated Barles--Souganidis framework to prove convergence to the unique viscosity solution (Theorem~\ref{thm:convergence}). For solutions with enhanced regularity, we further established an optimal-order error estimate of $O(\Delta t + h^2)$ (Theorem~\ref{thm:convergence_rate_narrative}).

Building upon this rigorous foundation, our second major contribution demonstrates the profound relevance of this numerical framework to a grand challenge in scientific machine learning. We proposed a novel Physics-Constrained Neural Operator (PCNO) architecture, whose two-stage design directly emulates the structure of our provably stable operator-splitting scheme. We then proved that this principled architecture is capable of breaking the curse of dimensionality for the difficult class of domain-to-solution mapping problems with changing topology (Theorem~\ref{thm:pcno_convergence_revised}). This result shows how the mathematical guarantees of a classical numerical method can provide an essential blueprint for constructing efficient and reliable neural operators.

The heart of this work is that the principled design of numerical schemes and the data-driven power of operator learning are not disparate pursuits but are deeply synergistic. By embedding the stable, convergent structure of a proven method into a learning architecture, we can overcome challenges that are intractable for purely data-driven or classical approaches alone.

This research opens several promising avenues for future investigation. From a numerical analysis perspective, natural extensions include the development of higher-order temporal discretizations to relax the CFL-type condition and the incorporation of a posteriori error estimators to drive adaptive mesh refinement. From a machine learning perspective, the immediate next step is the practical implementation and training of the proposed PCNO to numerically validate its theoretical efficiency on the perforated domain problem. Finally, extending the coupled analysis to a broader class of nonlinearities in the constraint operators remains an important and challenging open problem at the intersection of these fields.

\section*{Acknowledgement}
I acknowledge the financial support from the National Science and Technology Council of Taiwan (NSTC 111-2221-E-002-053-MY3).

\bibliography{sn-bibliography}% common bib file

\end{document}